\documentclass{amsart}

\usepackage{amssymb, amsmath}
\usepackage{mathrsfs}
\usepackage{amscd}

\usepackage[colorlinks,linkcolor={blue},
citecolor={blue},urlcolor={red},]{hyperref}

\theoremstyle{plain}
\newtheorem{theorem}{Theorem}[section]
\theoremstyle{remark}
\newtheorem{remark}[theorem]{Remark}
\newtheorem{example}[theorem]{Example}
\theoremstyle{plain}
\newtheorem{corollary}[theorem]{Corollary}
\newtheorem{lemma}[theorem]{Lemma}
\newtheorem{proposition}[theorem]{Proposition}

\numberwithin{equation}{section}

\def\Z{{\mathbb Z}}

\def\R{{\mathbb R}}

\newcommand{\E}{{\mathbb E}}
\renewcommand{\P}{{\mathbb P}}
\newcommand{\F}{{\mathcal F}}

\renewcommand{\a}{\alpha}
\renewcommand{\b}{\beta}
\newcommand{\g}{\gamma}

\renewcommand{\l}{\lambda}

\renewcommand{\O}{\Omega}

\newcommand{\beq}{\begin{equation}}
\newcommand{\eeq}{\end{equation}}
\newcommand{\bal}{\begin{aligned}}
\newcommand{\eal}{\end{aligned}}
\newcommand{\ben}{\begin{enumerate}}
\newcommand{\een}{\end{enumerate}}
\newcommand{\bit}{\begin{itemize}}
\newcommand{\eit}{\end{itemize}}

\newcommand{\bth}{\begin{theorem}}
\renewcommand{\eth}{\end{theorem}}
\newcommand{\bpr}{\begin{proposition}}
\newcommand{\epr}{\end{proposition}}
\newcommand{\ble}{\begin{lemma}}
\newcommand{\ele}{\end{lemma}}
\newcommand{\bpf}{\begin{proof}}
\newcommand{\epf}{\end{proof}}
\newcommand{\bex}{\begin{example}}
\newcommand{\eex}{\end{example}}
\newcommand{\bre}{\begin{example}}
\newcommand{\ere}{\end{example}}

\newcommand{\calL}{{\mathcal L}}
\newcommand{\n}{\Vert}
\newcommand{\one}{{{\bf 1}}}
\newcommand{\embed}{\hookrightarrow}
\newcommand{\s}{^*}
\newcommand{\lb}{\langle}
\newcommand{\rb}{\rangle}

\newcommand{\limn}{\lim_{n\to\infty}}

\newcommand{\sumk}{\sum_{k\ge 1}}
\newcommand{\sumj}{\sum_{j\ge 1}}

\newcommand{\wh}{\widehat}
\newcommand{\supp}{\text{\rm supp\,}}

\begin{document}

\title[Vector-valued Besov spaces and $\g$-radonifying operators]
{Embedding vector-valued Besov spaces into spaces of $\g$-radonifying
operators}

\author{Nigel Kalton}
\address{Mathematics Department\\
University of Missouri, Columbia \\ MO 65211 USA
\\ nigel@math.missouri.edu\\}

\author{Jan van Neerven}
\author{Mark Veraar}
\address{Delft Institute of Applied Mathematics\\
Technical University of Delft \\ P.O. Box 5031\\ 2600 GA Delft\\The
Netherlands \\ J.M.A.M.vanNeerven@tudelft.nl, M.C.Veraar@tudelft.nl\\}

\author{Lutz Weis}
\address{Mathematisches\, Institut\, I \\
\, Technische\, Universit\"at \, Karlsruhe \\
\, D-76128 \, Karls\-ruhe\\Germany \\
Lutz.Weis@mathematik.uni-karlsruhe.de}

\thanks{The second and third named
authors are supported by the `VIDI subsidie' 639.032.201 of the Netherlands
Organization for Scientific Research (NWO) and by the Research Training
Network HPRN-CT-2002-00281. The fourth named author was supported by grants
from the Volkswagenstiftung (I/78593) and the Deutsche Forschungsgemeinschaft
(We 2847/1-1)}



\begin{abstract}
It is shown that a Banach space $E$ has type $p$ if and only for some (all)
$d\ge 1$ the Besov space $B_{p,p}^{(\frac1p-\frac12)d}(\R^d;E)$ embeds into
the space $\g(L^2(\R^d),E)$ of $\g$-radonifying operators $L^2(\R^d)\to E$. A
similar result characterizing cotype $q$ is obtained. These results may be
viewed as $E$-valued extensions of the classical Sobolev embedding theorems.
\end{abstract}

\date\today

\maketitle

\section{Introduction}
Let $E$ be a real or complex Banach space and denote by ${\mathcal S}(\R^d;E)$
the Schwartz space of smooth, rapidly decreasing functions $f:\R^d\to E$. For
a function $f\in {\mathcal S}(\R^d;E)$ we consider the linear mapping $I_f:
L^2(\R^d)\to E$ defined by
$$ I_f g = \int_{\R^d} f(x)g(x)\,dx.$$
The aim of this paper is to prove the following characterization of
Banach spaces $E$ with type $p$ in terms of the embeddability of
certain $E$-valued Besov spaces
into spaces of $\g$-radonifying operators with values in $E$ and vice versa.
The precise definitions of the spaces $B_{p,p}^{(\frac{1}p-\frac{1}2)d}(\R^d;E)$ and $\g(L^2(\R^d),E)$
are recalled below.

\begin{theorem}\label{thm:main}
Let $E$ be a Banach space and let $1\le p\le 2\le q\le \infty$.
\ben
\item $E$ has type $p$ if and only if for
some (all) $d\ge 1$ the mapping $I: f\mapsto I_f$ extends to a continuous embedding
$$B_{p,p}^{(\frac{1}p-\frac{1}2)d}(\R^d;E)\embed  \g(L^2(\R^d),E);$$

\item $E$ has cotype $q$ if and only if for
some (all) $d\ge 1$ the mapping $I^{-1}: I_f\mapsto f$ extends to a continuous
embedding
$$\g(L^2(\R^d),E)\embed B_{q,q}^{(\frac{1}q-\frac{1}2)d}(\R^d;E).$$
\een
\end{theorem}

A version of this result for bounded open domains in $\R^d$ is obtained as well.

As is well known \cite{HJ-P,RS}, see also \cite{NW2}, $E$ has type $2$ if and
only if the mapping $f\mapsto I_f$ extends to a continuous embedding $
L^2(\R^d;E)\embed \g(L^2(\R^d),E)$, and $E$ has cotype $2$ if and only if $
\g(L^2(\R^d),E) \embed L^2(\R^d;E)$. Thus in some sense, Theorem
\ref{thm:main} may be viewed as an extension of these results for general
values of $p$ and $q$.

If $\dim E =1$, then
$\g(L^2(\R^d);E) = L^2(\R^d)$  and the embeddings
of Theorem \ref{thm:main}
reduce to the well-known Sobolev embeddings
$$ B_{p,p}^{(\frac{1}p-\frac{1}2)d}(\R^d)\embed L^2(\R^d)\embed
B_{q,q}^{(\frac{1}q-\frac{1}2)d}(\R^d), \quad 1\le p\le 2\le q\le \infty.$$

Vector-valued Besov spaces have attracted recent
attention in the theory of parabolic evolution equations in Banach spaces as a tool for
establishing optimal regularity results; see for instance \cite{Am,DHP}.
In \cite{GW1},
Fourier multiplier theorems with optimal exponents are established for operator-valued
multipliers on Besov spaces of functions taking values in Banach spaces with Fourier type $p$.

On the other hand, the spaces $\g(L^2(\R^d),E)$ have recently played an
important role in the theory of $H^\infty$-functional calculus for sectorial
operators \cite{FW, KW-Eucl, KW} and the theory of wavelet decompositions
\cite{KaW}. Furthermore, the spaces $\g(L^2(\R^d),E)$ have been characterized
in terms of stochastic integrals with respect to (cylindrical) Brownian
motions \cite{NVW, NW, RS}. Therefore, our results allow to compare various
square functions and they also give conditions for the stochastic
integrability of $E$-valued functions. These applications, which motivated our
results, will be detailed in a forthcoming paper.

\medskip
Throughout this paper, $H$ is a Hilbert space and $E$ is a Banach space,
which may be taken both real or both complex. Furthermore, $(r_n)_{n\ge 1}$ denotes a Rademacher sequence and $(\g_n)_{n\ge 1}$ a Gaussian sequence.

\subsection{Type and cotype}
Let $p\in [1,2]$ and $q\in [2,\infty]$.
A Banach space $E$ is said to have {\em type $p$} if there exists
a constant $C\ge 0$ such that for all finite subsets $\{x_1,\dots,x_N\}$ of $E$
we have
$$
\Big(\E \Big\n \sum_{n=1}^N r_n x_n\Big\n^2\Big)^\frac12
\le C \Big(\sum_{n=1}^N \n x_n\n^p\Big)^\frac1p.
$$
The least possible constant $C$ is called the {\em type $p$
constant} of $E$ and is denoted by $T_p(E)$.
A Banach space $E$ is said to have {\em cotype $q$} if there exists
a constant $C\ge 0$ such that for all finite subsets $\{x_1,\dots,x_N\}$ of $E$
we have
$$
\Big(\sum_{n=1}^N \n x_n\n^q\Big)^\frac1q \le C\Big(\E \Big\n \sum_{n=1}^N r_n x_n\Big\n^2\Big)^\frac12,
$$
with the obvious modification in the case $q=\infty$. The least possible
constant $C$ is called the {\em cotype $q$ constant} of $E$ and is denoted by
$C_p(E)$. As is well known, in both definitions the r\^ole of the Rademacher
variables may be replaced by Gaussian variables without altering the class of
spaces under consideration. The least constants arising from these equivalent
definitions are called the {\em Gaussian type $p$ constant} and the {\em
Gaussian cotype $q$ constant} of $E$ respectively, notation $T_p^\g(E)$ and
$C_q^\g(E)$.

Every Banach space has type $1$ and cotype $\infty$. The $L^p$-spaces have
type $\min\{p,2\}$ and cotype $\max\{p,2\}$ for $1\le p<\infty$.
Every Hilbert space has both type $2$ and cotype $2$, and a famous result of Kwapie\'n asserts that up to isomorphism this property characterizes the class of Hilbert spaces.

For more information we refer to Maurey's survey article \cite{Mau} and the references given therein.

\subsection{Besov spaces}\label{subsec:besov}
Next we recall the definition of Besov spaces using the so-called Little\-wood-Paley decomposition. We follow the approach of Peetre; see \cite[Section 2.3.2]{Tr} (where the scalar-valued case is considered) and \cite{Am,GW1,Schm}.
The Fourier transform of a function $f\in L^1(\R^d;E)$ will be
normalized  as
$$ \wh{f}(\xi) = \frac1{(2\pi)^{d/2}}\int_{\R^d} f(x)e^{-ix\cdot\xi}\,dx, \quad \xi\in\R^d.$$

Let $\phi\in{\mathscr S}(\R^d)$ be a fixed Schwartz function
whose Fourier transform $\wh\phi$
is nonnegative and has support in $\{\xi\in\R^d: \ \tfrac12\le |\xi|\le 2\}$ and which satisfies
$$ \sum_{k\in\Z} \wh\phi(2^{-k}\xi) =1 \quad\hbox{for $\xi\in \R^d\setminus\{0\}$}.$$
Define the sequence $({\varphi_k})_{k\ge 0}$ in ${\mathscr S}(\R^d)$ by
$$\wh{\varphi_k}(\xi) = \wh\phi(2^{-k}\xi) \quad \text{for}\ \  k=1,2,\dots \quad \text{and} \ \ \wh{\varphi_0}(\xi) = 1- \sum_{k\ge 1} \wh{\varphi_k}(\xi), \quad \xi\in\R^d.$$

For $1\le p,q \le \infty$ and $s\in\R$ the {\em Besov space} $B_{p,q}^s(\R^d;E)$ is defined as the space of all $E$-valued tempered distributions
$f\in {\mathscr S}'(\R^d;E)$ for which
$$ \n f\n_{B_{p,q}^s (\R^d;E)} := \Big\n \big( 2^{ks}{\varphi}_k * f\big)_{k\ge 0} \Big\n_{l^q(L^p(\R^d;E))} $$
is finite.
Endowed with this norm, $B_{p,q}^s(\R^d;E)$ is a Banach space, and up to an equivalent norm this space is independent of the choice of the initial function  $\phi$.
The sequence $({\varphi}_k * f)_{k\ge 0}$ is
called the {\em Littlewood-Paley decomposition} of $f$ associated with the function $\phi$.

The following continuous inclusions hold:
\[\ B_{p, q_1}^s(\R^d; E) \embed B_{p, q_2}^s(\R^d;E), \ B_{p, q}^{s_1}(\R^d;
E) \embed B_{p, q}^{s_2}(\R^d;E)\] for all $s,s_1, s_2\in \R$, $p,
q, q_1, q_2\in [1, \infty]$ with $q_1\leq q_2$, $s_2 \leq s_1$.
Also note that
$$ B_{p, 1}^0(\R^d; E)\embed L^p(\R^d;E)\embed B_{p, \infty}^0(\R^d; E).$$
If $1\le p,q<\infty$, then $B_{p,q}^s(\R^d;E)$ contains the Schwartz space ${\mathscr S}(\R^d;E)$ as a dense subspace.

In Section \ref{sec:besovD} we shall need the following lemma.
For $\l>0$ let $f_\l(x) := f(\l x)$.

\begin{lemma}\label{lem:dilation}
Let $p,q\in [1,\infty]$ and $s\in\R$, $s\not=0$.
\ben
\item
If $s > 0$, there exists a constant $C>0$
such that for all
$\l = 2^n$, $n\ge 1$, and $ f\in B_{p,q}^s(\R^d;E)$ we have
$$  \n f_\l \n_{B_{p,q}^s(\R^d;E)} \le C\l^{s-\frac{d}{p}}
\n f \n_{B_{p,q}^s(\R^d;E)}.
$$
\item
If $s< 0$, there exists a constant $C>0$
such that for all
$\l = 2^n$, $n\le -1$, and $ f\in B_{p,q}^s(\R^d;E)$ we have
$$  \n f_\l \n_{B_{p,q}^s(\R^d;E)} \le C\l^{s-\frac{d}{p}}
\n f \n_{B_{p,q}^s(\R^d;E)}.
$$
\een
\end{lemma}
\begin{proof}
We only prove (1), the proof of (2) being similar. The proofs are
patterned after \cite[Proposition 3.4.1]{Tr2}.

Let $\phi$ and $\varphi_k$, $k=0,1,2,\dots,$ be as in Subsection
\ref{subsec:besov}. Define, for $m\in\Z$, the functions $\psi_m$ by
$\wh{\psi_m}(\xi) := \wh\phi(2^{-m}\xi)$. Then $\psi_m = \varphi_m$ for
$m=1,2,\dots$ and $(\wh{\psi_m})_\l=\wh{\psi_{m-n}}$ for $m\in\Z$ and
$\l=2^n$, $n\in\Z$. For $s> 0$ we have
$$
\bal
  \Big(\sum_{k\ge 0} 2^{ksq} \n \varphi_k * f_\l\n_{L^p(\R^d;E)}^q\Big)^\frac1q
& =  \l^{-\frac{d}{p}}\Big(\sum_{k\ge 0} 2^{ksq} \n \F^{-1}((\wh{\varphi_k})_\l \wh f)\n_{L^p(\R^d;E)}^q\Big)^\frac1q
\\ & \le
 \l^{-\frac{d}{p}} \n \F^{-1}((\wh{\varphi_0})_\l \wh f)\n_{L^p(\R^d;E)}
\\ & \qquad +
 \l^{-\frac{d}{p}}\Big(\sum_{k=1}^n 2^{ksq} \n \F^{-1}(\wh{\psi_{k-n}} \wh f)\n_{L^p(\R^d;E)}^q\Big)^\frac1q
\\ & \qquad + \l^{s-\frac{d}{p}}\Big(\sum_{l\ge 1} 2^{lsq}
\n \F^{-1}(\wh{\psi_l}\wh f)\n_{L^p(\R^d;E)}^q\Big)^\frac1q
\\ &  =: {\rm (I)} + {\rm (II)}+ {\rm (III)}.
\eal
$$
Since $\wh{\varphi_0}=1$ on $(0,1]$ and $(\wh{\varphi_0})_\l$ has support in
$(0,2^{-n}]\subseteq (0,\frac12]$, by Young's inequality we have
$$
\n \F^{-1}((\wh{\varphi_0})_\l \wh f)\n_{L^p(\R^d;E)}
 = \n  \F^{-1}((\wh{\varphi_0})_\l \wh{\varphi_0}\wh f)\n_{L^p(\R^d;E)}
\le  \n \varphi_0\n_{L^1(\R^d)} \n \varphi_0 * f\n_{L^p(\R^d;E)}.
$$
Hence, $${\rm (I)} \le \l^{-\frac{d}{p}} \n \varphi_0\n_{L^1(\R^d)} \n f \n_{B_{p,p}^s(\R^d;E)}
 \le \l^{s-\frac{d}{p}} \n \varphi_0\n_{L^1(\R^d)} \n f \n_{B_{p,q}^s(\R^d;E)}.$$
To estimate (II) we note that for $k=1,\dots, n-1$ the functions
$\wh{\psi_{k-n}}$ have support in $(0,1]$. Therefore,
$$\bal
 \n \F^{-1}(\wh{\psi_{k-n} }\wh f)\n_{L^p(\R^d;E)}
 &  \le \n \psi_{k-n}\n_{L^1(\R^d)} \n \varphi_0 * f\n_{L^p(\R^d;E)}
\\ &  = \n \phi\n_{L^1(\R^d)} \n \varphi_0 * f\n_{L^p(\R^d;E)}.
\eal $$ Similarly, for $k=n$,
$$
\bal \ & \n \F^{-1}((\wh{\psi_n})_\l \wh f)\n_{L^p(\R^d;E)}
  \le  \n \phi\n_{L^1(\R^d)} \big(\n \varphi_0 * f\n_{L^p(\R^d;E)}
+\n \varphi_1 * f\n_{L^p(\R^d;E)} \big). \eal
$$
Summing these terms and using that $s> 0$ we obtain
$${\rm (II)} \le  C_{s,q}\l^{s-\frac{d}{p}}\n \varphi\n_{L^1(\R^d)}
\n f \n_{B_{p,q}^s(\R^d;E)}$$
with a constant $C_{q,s}$ depending only of $q$ and $s$.
Obviously, $$ {\rm (III)} \le \l^{s-\frac{d}{p}} \n f \n_{B_{p,q}^s(\R^d;E)}.$$
By putting these estimates together the desired inequality follows.
\end{proof}

\subsection{$\g$-Radonifying operators}
For a finite rank operator $R:H\to E$ of the form
\beq\label{eq:T} Rh = \sum_{n=1}^N [h,h_n]_H\,x_n \eeq
with $h_1,\dots,h_N$ orthonormal in $H$, we define
$$ \n R\n_{\g(H,E)}^2 := \E \Big\n \sum_{n=1}^N \g_n Rh_n\Big\n^2.$$
Note that $\n R\n_{\g(H,E)}$ does not depend on the particular representation of $R$ as in \eqref{eq:T}.
The completion of the space of finite rank operators with respect to the norm $\n\cdot\n_{\g(H,E)}$ defines a two-sided operator ideal $\g(H,E)$ in $\calL(H,E)$. If $H$ is separable, an operator $R\in\calL(H,E)$
belongs to $\g(H,E)$ if and only if for some (equivalently,
for every) orthonormal basis $(h_n)_{n\ge 1}$ of $H$ the Gaussian sum
$\sum_{n\ge 1} \g_n Rh_n$ converges in $L^2(\O;E)$, in which case we have
$$\n R\n_{\g(H,E)}^2 = \E\Big\n \sum_{n\ge 1} \g_n Rh_n\Big\n^2.$$
We refer to \cite[Chapter 12]{DJT} for more information.

The following elementary convergence result, cf. \cite[Proposition 2.4]{NVW}, will be useful.
If the $T_1, T_2, \ldots\in \calL(H)$ and $T\in \calL(H)$ satisfy
$\sup_{n\ge 1}\|T_n\| < \infty$ and $\limn T^* h = T_n^* h$ for all $h\in H$,
then for all $R\in \g(H,E)$ we have
\beq\label{eq:conv-g}
\limn \n R \circ T_n - R \circ T\n_{\g(H,E)}=0.
\eeq

If $H_1$ and $H_2$ are Hilbert spaces, then every bounded operator $T: H_1\to H_2$
induces a bounded operator $\widetilde T: \g(H_1,E)\to \g(H_2,E)$ by the formula
$$\widetilde T R := R\circ T\s$$
and we have
\begin{equation}\label{KW}
\n \widetilde T\n_{\calL(\g(H_1,E),\g(H_2,E))}  \le \n T\n_{\calL(H_1,H_2)}.\end{equation}
This extension procedure is introduced in \cite{KW} and will be
useful below.

If $(S,\Sigma,\mu)$ is a $\sigma$-finite measure space, we denote by
$\g(S;E)$ the vector space of all strongly $\mu$-measurable functions
$f:S\to E$ for which $\lb f,x\s\rb$ belongs to $L^2(S)$ for all $x\s\in E\s$ and the associated Pettis operator $I_f: L^2(S)\to E$,
$$ I_f g = \int_S fg\,d\mu$$
belongs to $\g(L^2(S),E)$. We identify functions defining the same operator.
An easy approximation argument shows that the simple functions in $\g(S;E)$
form a dense subspace of $\g(L^2(S),E)$. We shall write
$$\n f\n_{\g(S;E)} := \n I_f\n_{\g(L^2(S),E)}.$$

\section{Embedding results for $\R^d$}\label{sec:besov}

The proof of Theorem \ref{thm:main} is based on two lemmas.

\begin{lemma}\label{lem:LN1} \
\ben
\item Let $E$ have type $p\in [1,2]$. If $f\in {\mathscr S}(\R^d;E)$ satisfies $\supp \wh f \subseteq [-\pi,\pi]^d$,
then $f\in \g(\R^d;E)$ and
$$ \n f\n_{\g(\R^d;E)} \le T_p^\g(E)  \n f\n_{L^p(\R^d;E)},
$$
where $T_p^\g(E)$ denotes the Gaussian type $p$ constant of $E$.
\item Let $E$ have cotype $q\in [2,\infty]$. If $f\in {\mathscr S}(\R^d;E)$ satisfies $\supp \wh f \subseteq [-\pi,\pi]^d$, then
$$
\n f\n_{\g(\R^d;E)} \ge C_q^\g(E)^{-1} \n f\n_{L^q(\R^d;E)},
$$
where $C_q^\g(E)$ denotes the Gaussian cotype $q$ constant of $E$.
\een
\end{lemma}
\begin{proof}
Let $Q := [-\pi,\pi]^d$.
We consider the functions $h_n(x) = (2\pi)^{-d/2} e^{in\cdot x}$
with $n\in\Z^d$, $x\in Q$, which define an orthonormal basis for $L^2(Q)$.

(1) \
Define the bounded operators $I_{f}: L^2(\R^d)\to E$
and $I_{\wh f}: L^2(\R^d)\to E$ by
$$ I_{f} g := \int_{\R^d} f(x) g(x)\,dx, \quad I_{\wh f} g := \int_{\R^d} \wh f(x) g(x)\,dx.$$
In case $E$ is a real Banach space
we consider its complexification in the second definition.
By the assumption on the support of $\wh f$ we may identify $I_{\wh{f}}$ with a bounded operator from $L^2(Q)$ to $E$ of the same norm.
Since $I_{\wh{f}} h_n = f(n)$,
for any finite subset $F\subseteq \Z^d$ we have
$$
\bal
\Big(\E \Big\n \sum_{n\in F} \g_n I_{\wh{f}}h_n \Big\n^2\Big)^\frac12
 = \Big(\E \Big\n \sum_{n\in F} \g_n f(n)\Big\n^2\Big)^\frac12
\le T_p^\g(E) \Big(\sum_{n\in F} \n f(n)\n^p\Big)^\frac1p.
\eal
$$
It follows that $I_{\wh{f}}\in \g(L^2(Q),E)$. By
the identification made above it follows that $I_{\wh{f}}\in \g(L^2(\R),E)$ and
$$
 \n I_{\wh f}\n_{\g(L^2(\R^d),E)}
 = \n I_{\wh f}\n_{\g(L^2(Q),E)}
\le  T_p^\g(E) \Big(\sum_{n\in \Z^d} \n f(n)\n^p\Big)^\frac1p.
$$
From \eqref{KW} it follows that
$$
\n f\n_{\g(\R^d;E)} =
\n I_f\n_{\g(L^2(\R^d),E)}
 = \n I_{\wh f}\n_{\g(L^2(\R^d),E)}
\le  T_p^\g(E) \Big(\sum_{n\in \Z^d} \n f(n)\n^p\Big)^\frac1p.
$$
For $t\in R := [0,1]^d$ put $f_t(s) = f(s+t)$.
Then $\supp \wh {f_t} \subseteq Q$ and
$$ \n f\n_{\g(\R^d;E)} = \n f_t\n_{\g(\R^d;E)}
\le T_p^\g(E) \Big(\sum_{n\in \Z^d} \n f_t(n)\n^p\Big)^\frac1p.
$$
By raising both sides to the power $p$ and integrating over $R$ we obtain
$$
\n f\n_{\g(\R^d;E)}
 \le T_p^\g(E)\Big( \int_R \sum_{n\in \Z^d} \n f_t(n)\n^p\,dt\Big)^\frac1p
= T_p^\g(E)\Big( \int_{\R^d} \n f(s)\n^p\,ds\Big)^\frac1p.
$$

(2) \ This is proved similarly. Note that by part (1) (with $p=1$) we have
$f\in\g(\R^d;E)$.
\end{proof}

Let $(S,\Sigma,\mu)$ be a measure space. For a bounded operator $R:L^2(S)\to
E$ and a set $S_0\in\Sigma$ we define $R|_{S_0}: L^2(S)\to E$ by $$ R|_{S_0}g
:= R(\one_{S_0}g).$$ Note that if $R\in\g(L^2(S),E)$, then
$R|_{S_0}\in\g(L^2(S),E)$ and
$$\n R|_{S_0}\n_{\g(L^2(S),E)} \le \n R\n_{\g(L^2(S),E)}$$
by the operator ideal property of $\g(L^2(S),E)$.

In the following lemma we use the well known fact that if $E$ has type $p$ (cotype $q$), then the same is true for the space $L^2(\O;E)$
and we have $$ T_p(L^2(\O;E)) = T_p(E), \quad  C_q(L^2(\O;E)) = C_q(E).$$

\begin{lemma} \label{lem:LN2}
Let $(S,\Sigma,\mu)$ be a measure space and let $(S_j)_{j\ge 1}\subseteq
\Sigma$ be a partition of $S$. \ben \item Let $E$ have type $p\in [1,2]$. Then
for all $R\in \g(L^2(S),E)$ we have
$$ \n R\n_{\g(L^2(S),E)} \le T_p(E)\Big(\sumj  \n R|_{S_j}\n_{\g(L^2(S),E)}^p\Big)^\frac1p.$$
\item Let $E$ have cotype $q\in [2,\infty]$.
Then for all $R\in \g(L^2(S),E)$ we have
$$  \n R\n_{\g(L^2(S),E)} \ge C_q(E)^{-1} \Big(\sumj  \n R|_{S_j}\n_{\g(L^2(S),E)}^q\Big)^\frac1q.$$
\een
\end{lemma}

\begin{proof}
(1) \
We may assume that $\mu(S_j)>0$ for all $j$.
Fixing $R$, we may also assume that  $\Sigma$ is countably generated.
As a result, $L^2(S)$ is separable and we may
choose an orthonormal basis $(h_{jk})_{j,k\ge 1}$ for $L^2(S)$ in such a way
that for each $j$ the sequence $(h_{jk})_{k\ge 1}$ is an orthonormal basis for
$L^2(S_j)$.  Let $(\g_{jk})_{j,k\ge 1}$ and $(r_j')_{j\ge 1}$ be a doubly-indexed Gaussian sequence and a  Rademacher sequence on probability spaces
$(\O,\P)$ and $(\O',\P')$, respectively. By a standard randomization argument,
$$\bal \n R\n_{\g(L^2(S),E)}
& = \Big(\E \Big\n \sum_{j,k\ge 1} \g_{jk} R h_{jk}
\Big\n^2\Big)^\frac12
\\ & = \Big(\E \Big\n \sum_{j,k\ge 1} \g_{jk} R|_{S_j} h_{jk}
\Big\n^2\Big)^\frac12
\\ & =\Big(\E' \Big\n \sum_{j\ge 1} r_j' \sumk  \g_{jk}R|_{S_j} h_{jk}
\Big\n_{L^2(\O;E)}^2\Big)^\frac12
\\ & \le T_p(L^2(\O;E))
\Big(\sum_{j\ge 1} \Big\n\sumk  \g_{jk}R|_{S_j} h_{jk} \Big\n_{L^2(\O;E)}^p\Big)^\frac1p
\\ & =  T_p(E) \Big(\sumj  \n R|_{S_j}
\n_{\g(L^2(S),E)}^p\Big)^\frac1p.
\eal
$$

(2) \ This is proved similarly.
\end{proof}

We are now prepared for the proof of Theorem \ref{thm:main}.
Recall that the Schwartz functions $\phi$ and $\varphi_k$, $k\ge 1$, are
defined in Subsection \ref{subsec:besov}.

\begin{proof}[Proof of Theorem \ref{thm:main}]
(1) \ First we prove the `only if' part and assume that $E$ has type $p$.
Let $f\in {\mathscr S}(\R^d;E)$ and let
$f_k := \varphi_k*f$.
Putting $g_k(x) := f_k(2^{-k}x)$ we have $g_k\in{\mathscr S}(\R^d;E)$ and
$$\supp \wh{g_k} \subseteq \{\xi\in \R^d: \ \tfrac12 \le |\xi| \le 2\}\subseteq [-\pi,\pi]^d.$$ Hence from Lemma \ref{lem:LN1} we obtain $f_k\in \g(\R^d;E)$ and
$$
\bal
\n f_k \n_{\g(\R^d;E)}
& = 2^{{-kd}/{2}} \n g_k \n_{\g(\R^d;E)}
\\ & \le 2^{{-kd}/{2}}T_p^\g(E) \n g_k \n_{L^p(\R^d;E)}
= 2^{\frac{kd}{p}-\frac{kd}{2}} T_p^\g(E) \n f_k \n_{L^p(\R^d;E)}.
\eal
$$
Using the Lemma \ref{lem:LN2}, applied to the decompositions $(S_{2k})_{k\in \Z}$
and $(S_{2k+1})_{k\in \Z}$ of $\R^d\setminus\{0\}$, we obtain, for all
$n\ge m\ge 0$,
$$
\bal
 \Big\n \sum_{k=2m}^{2n} f_k \Big\n_{\g(\R^d;E)}
 & \le  T_p^\g(E)T_p(E) \Big(\sum_{j= m}^n 2^{(\frac{2jd}{p} - \frac{2jd}{2})p} \n {f_{2j}}\n_{L^p(\R^d;E)}^p\Big)^\frac1p
\\ & \qquad +
T_p^\g(E)T_p(E)  \Big(\sum_{j=m}^{n-1} 2^{(\frac{(2j+1)d}{p} - \frac{(2j+1)d}{2})p} \n {f_{2j+1}}\n_{L^p(\R^d;E)}^p  \Big)^\frac1p.
\eal
$$
Estimating sums of the form $\sum_{k=2m}^{2n+1}$, $ \sum_{k=2m+1}^{2n}$,  and $\sum_{k=2m+1}^{2n+1}$ in a similar way,
it follows that $f\in \g(\R^d;E)$ and
$$ \n f\n_{\g(\R^d;E)} \le 2T_p^\g(E)T_p(E) \n f\n_{B_{p,p}^{(\frac{1}{p} - \frac{1}{2})d}(\R^d;E)}.
$$
Since ${\mathscr S}(\R^d;E)$ is dense in $B_{p,p}^{(\frac{1}{p} - \frac{1}{2})d}(\R^d;E)$ it follows that the mapping $f\mapsto I_f$ extends to a bounded operator $I$ from $B_{p,p}^{(\frac{1}{p} - \frac{1}{2})d}(\R^d;E)$
into $\g(\R^d;E)$ of norm
$\n I\n \le 2T_p^\g(E) T_p(E)$.
The simple proof that $I$ is injective is left to the reader.

Next we prove the `if' part. For $n\geq 1$, let $\psi_n\in \mathcal{S}(\R^d)$
be defined as
\[\widehat{\psi}_n(\xi) = c 2^{-nd/2}\widehat{\phi}(2^{-n}\xi),\]
where $c:=\|\phi\|_{L^2(\R^d)}^{-1}$. Then $(\psi_{3n})_{n\geq 1}$ is an
orthonormal system in $L^2(\R^d)$. For any finite sequence $(x_n)_{n=1}^N$ in
$E$ we then have, with $f:=\sum_{n=1}^N \psi_{3n} \otimes x_n$,
\[\|f\|_{\g(\R^d;E)}^2 = \E\Big\|\sum_{n=1}^N \g_n x_n\Big\|^2.\]
Notice that for $k\geq 1$,
\[\|\varphi_{k}*\varphi_{k}\|_{L^p(\R^d)} = 2^{kd - \frac1p k d}  \|\phi*\phi\|_{L^p(\R^d)} \]
and
\[\|\varphi_{k+1}*\varphi_{k}\|_{L^p(\R^d)} = 2^{kd - \frac1p k d}  \|\varphi_1*\phi\|_{L^p(\R^d)}.\]
Therefore, for $n=1, \ldots, N$,
\[\|\varphi_{3n}*f\|_{L^p(\R^d;E)} = c 2^{-\frac32 nd} \|\varphi_{3n}*\varphi_{3n}\|_{L^p(\R^d)} \|x_n\| = c 2^{(\frac12-\frac1p) 3nd}  \|\phi*\phi\|_{L^p(\R^d)} \|x_n\| \]
and similarly,
\[\begin{aligned}
\|\varphi_{3n-1}*f\|_{L^p(\R^d;E)}
= c2^{(\frac12-\frac1p)3nd  -(1-\frac1p)d}
\|\varphi_1*\phi\|_{L^p(\R^d)} \|x_n\|
\end{aligned}\]
and
\[\begin{aligned}
\|\varphi_{3n+1}*f\|_{L^p(\R^d;E)}
= c2^{(\frac12-\frac1p)3nd} \|\varphi_1*\phi\|_{L^p(\R^d)} \|x_n\|.
\end{aligned}\]
Finally, for $k\ge 3N+2$ we have $\varphi_{k}*f = 0$.  Summing up, it follows that there
exists a constant $C$, depending only on $p,d$ and $\phi$ such that
\[\begin{aligned}
\|f\|_{B_{p,p}^{(\frac1p-\frac12)d}(\R^d;E)} \leq C \Big(\sum_{n = 1}^N
\|x_n\|^p\Big)^{\frac1p}.
%
\end{aligned}\]
By putting things together we see that $E$ has type $p$, with Gaussian type
$p$ constant
$
T_p^\g(E)\le  C\n I\n,
$
where  $I: B_{p,p}^{(\frac{1}{p} -
\frac{1}{2})d}(\R^d;E)\embed \g(\R^d;E)$ is the embedding.


(2) This is proved similarly.
\end{proof}

As a special case of Theorem \ref{thm:main},
note that for every Banach space $E$ we obtain continuous embeddings
$$ B_{1,1}^{ \frac{1}{2}d}(\R^d;E) \embed \g(L^2(\R^d),E)\embed
 B_{\infty,\infty}^{ - \frac{1}{2}d}(\R^d;E).
$$
As is easily checked by going through the proofs, these embeddings are contractive.

Let $H^{\b,p}(\R^d;E)$, with $\b\in\R$ and $1\le p<\infty$, denote the usual $E$-valued Lebesgue-Bessel potential spaces \cite[Section 6.2]{BL}, \cite[Section 2.33]{Tr}.
In \cite{KW-Eucl} the $\g$-Sobolev spaces $\g(H^{\b,2}(\R^d),E)$
are introduced and their basic properties are studied. From Theorem \ref{thm:main} we obtain the following $\g$-analogue of the Sobolev embedding theorem.

\begin{corollary}\label{cor:mainRd} \
\ben
\item
If $E$ has type $p\in [1,2]$, we have continuous embeddings
$$ H^{\a,p}(\R^d;E) \embed B_{p,p}^{\b+(\frac{1}{p} - \frac{1}{2})d}(\R^d;E)\embed \g(H^{\b,2}(\R^d),E)$$
for all $\a,\b\in\R$ satisfying $\a> \b +(\frac{1}{p} - \frac{1}{2})d.$
\item If $E$ has cotype $q\in [2,\infty]$, we have continuous embeddings
$$\g(H^{\b,2}(\R^d),E) \embed B_{q,q}^{\b+(\frac{1}{q} - \frac{1}{2})d}(\R^d;E)
\embed H^{\a,q}(\R^d;E)$$ for all $\a, \b\in\R$ satisfying $\a < \b
+(\frac{1}{q} - \frac{1}{2})d.$ \een
\end{corollary}

\begin{remark} Taking $q=\infty$ in (2), as a special case we obtain
the embedding $\g(H^{\b,2}(\R^d),E) \embed
B_{\infty,\infty}^{\b-\frac{d}{2}}(\R^d;E).$ If $\b-\frac{d}{2}$ is strictly
positive and not an integer, the latter space can be identified, up to an
equivalent norm, with the H\"older space $(BUC)^{\b-\frac{d}{2}}(\R^d;E)$
\cite[Equation (5.8)]{Am} and we thus obtain a continuous embedding
\beq\label{eq:BUC} \g(H^{\b,2}(\R^d);E) \embed (BUC)^{\b-\frac{d}{2}}(\R^d;E).
\eeq
\end{remark}


\begin{proof} The second embedding in (1) and the first embedding in (2) are immediate from Theorem \ref{thm:main} combined with the fact that $(I-\Delta)^{-\b/2}$
acts as an isomorphism from $B_{q,q}^{(\frac{1}{q} - \frac{1}{2})d}(\R^d;E)$ onto $B_{q,q}^{\b+(\frac{1}{q} - \frac{1}{2})d}(\R^d;E)$ \cite[Theorem 6.1]{Am} and, by \eqref{KW},
from $\g(L^2(\R^d),E)$ onto $\g(H^{\b,2}(\R^d),E)$. The first embedding in (1) and the second embedding in (2) follow from the $E$-valued analogues of \cite[Theorem 6.2.4]{BL}. 
\end{proof}

Note that (2) can be combined with the classical Sobolev embedding theorem
to yield an inclusion result which is slightly weaker than \eqref{eq:BUC}.

\section{Embedding results for bounded domains}\label{sec:besovD}

Let $D$ be a nonempty bounded open domain in $\R^d$. For
$1\le p,q\le \infty$ and $s\in\R$ we define
$$B_{p,q}^s (D;E) = \{f|_D: \ f\in B_{p,q}^s (\R^d;E)\}.$$
This space is a Banach space endowed with the norm
$$ \n g\n_{  B_{p,q}^s (D;E)} = \inf_{f|_D = g} \n f\n_{B_{p,q}^s (\R^d;E)}.$$
See \cite[Section 3.2.2]{Tr2} (where the scalar case is considered) and \cite{Am2}.

In Theorem \ref{thm:embed} below we shall
obtain a version of Theorem \ref{thm:main} for bounded domains.
We need the following lemma, where for $r>0$ we denote $B_r := \{x\in E: \|x\|<r\}$.
\begin{lemma}\label{lem:compactsupp}
Let $1\le p,q\le \infty$, $s\in \R$. There exists a constant $C$ such that for
every $r\geq 1$ and for all $f\in B^{s}_{p,q}(\R^d;E)$ with
$\supp(f)\subseteq B_r$,
\[\|f\|_{B^{s}_{p,q}(\R^d;E)}\leq C \|f|_{B_{2r}}\|_{B^{s}_{p,q}(B_{2r};E)}.\]
\end{lemma}
\begin{proof}
Choose $\psi\in {\mathscr S}(\R^d)$ such that $\psi\equiv 1$ on $B_1$ and
$\psi\equiv 0$ outside $B_{2}$. Fix an integer
$k>\max\big\{s,\frac{d}{p}-s\big\}$. Notice that for the
$\frac{1}{r}$-dilation $\psi_{\frac1r}(x) := \psi(\frac1r x)$ we have
$\|\psi_{\frac{1}{r}}\|_{W^{k,\infty}(\R^d)}
\leq\|\psi\|_{W^{k,\infty}(\R^d)}$.
Choose $g\in B^{s}_{p,q}(\R^d;E)$ such that $g \equiv f$ on $B_{2r}$ and
\[\|g\|_{B^{s}_{p,q}(\R^d;E)}\leq 2\|f|_{B_{2r}}\|_{B^{s}_{p,q}(B_{2r};E)}.\]
Then it follows from the vector-valued generalization of \cite[Theorem
2.8.2]{Tr2} that
\[\begin{aligned}
\|f\|_{B^{s}_{p,q}(\R^d;E)}&=\|\psi_{\frac1r}f\|_{B^{s}_{p,q}(\R^d;E)}
=\|\psi_{\frac1r}g\|_{B^{s}_{p,q}(\R^d;E)}
\\ & \leq c\|\psi\|_{W^{k,\infty}(\R^d)} \|g\|_{B^{s}_{p,q}(\R^d;E)}
\leq C\|f|_{B_{2r}}\|_{B^{s}_{p,q}(B_{2r};E)},
\end{aligned}\]
where $C = 2 c\|\psi\|_{W^{k,\infty}(\R^d)}$.
\end{proof}


\begin{theorem} \label{thm:embed}
Let $1\leq p\leq 2\leq q\leq \infty$ and let $D\subseteq \R^d$ be a nonempty bounded open domain.
\begin{enumerate}
\item $E$ has type $p$ if and only if we have a
continuous embedding
$$ B_{p,p}^{(\frac{1}{p} - \frac{1}{2})d}(D;E) \embed \g(L^2(D),E).$$
\item $E$ has cotype $q$ if and only if we have a continuous embedding
$$ \g(L^2(D),E) \embed B_{q,q}^{(\frac{1}{q} - \frac{1}{2})d}(D;E).$$
\end{enumerate}
In both cases, the norm of the embedding does not exceed the norm of the
corresponding embedding with $D$ replaced by $\R^d$.
\end{theorem}
Note again the special cases corresponding to $p=1$ and $q=\infty$, which hold for arbitrary  Banach spaces $E$.
Corollary \ref{cor:mainRd} admits a version for bounded domains as well.
\begin{proof}
The ``only if'' parts in (1) and (2) and the final remark follow directly from
the definition.

For the proofs of the ``if'' parts in (1) and (2), there is no loss of generality in assuming that $0\in D$. Let $D_n = 2^nD$ and note that $\one_{D_n} \to
\one$ pointwise. The idea is to `dilate' the embedding for $D$ to $D_n$ and
pass to the limit $n\to\infty$ to obtain the corresponding embedding for
$\R^d$. That $E$ has type $p$ or cotype $q$ is then a consequence of Theorem
\ref{thm:main}.

(1): The result being trivial for $p=1$ we shall assume that $p\in (1,2]$. Fix a
function $f\in {\mathscr S}(\R^d;E)$ and note that by Lemma \ref{lem:LN1}
(applied with $p=1$) that $f\in \g(\R^d;E)$. Fix $n\ge 1$ arbitrary and put
$f_n := f|_{D_n}$. Let $g_n:D\to E$ be defined by
$$g_n(x):= f_n(2^nx), \qquad x\in D.$$
Then $g_n\in \g(D;E)$ and
$$ \n g_n\n_{\g(D;E)} = 2^{-\frac{1}{2}nd} \n f_n\n_{\g(D_n;E)}. $$
Also, $g_n = g^{(n)}|_D$, where $g^{(n)}(x) = f(2^nx)$ for $x\in\R^d$. By
Lemma \ref{lem:dilation} there exists a constant $C>0$, independent of $n$,
such that
$$ \n g_n\n_{ B_{p,p}^{(\frac1p-\frac12)d}(D;E)}
\le  \n g^{(n)}\n_{ B_{p,p}^{(\frac1p-\frac12)d}(\R^d;E)} \le
C2^{-\frac{1}{2}nd} \n f\n_{B_{p,p}^{(\frac1p-\frac12)d}(\R^d;E)}.$$
Denoting by $I: B_{p,p}^{(\frac1p-\frac12)d}(D;E)\embed \g(D;E)$  the
embedding, it follows that
$$
\bal \n f_n\n_{\g(D_n;E)} & = 2^{\frac12 nd}\n g_n\n_{\g(D;E)}
\\ & \le  2^{\frac12 nd}\n I\n \, \n g_n\n_{B_{p,p}^{(\frac1p-\frac12)d}(D;E)}
 \le C \n I\n \, \n f\n_{B_{p,p}^{(\frac1p-\frac12)d}(\R^d;E)}.
\eal
$$
Passing to the limit $n\to \infty$ we obtain, by virtue of \eqref{eq:conv-g},
$$ \n f\n_{\g(\R^d;E)}\le C \n I\n \, \n f\n_{B_{p,p}^{(\frac1p-\frac12)d}(\R^d;E)}.
$$
An application of Theorem \ref{thm:main} finishes the proof.


(2): It suffices to consider the case $q\in [2,\infty)$. Fix $f\in
C^{\infty}_c(\R^d;E)$ and let $r\geq 1$ be so large that $\supp(f) \subseteq
B_r$. With the same arguments as in (1) one can show that
\[\|f_n\|_{B_{q,q}^{(\frac1q-\frac12)d}(D_n;E)}\leq C \|f\|_{\g(\R^d;E)},\]
where $f_n = f|_{D_n}$ as before and $C$ is a constant not depending on $f$
and $n$. It follows from Lemma \ref{lem:compactsupp} that there is a constant
$C'$, independent of $f$ and $r$, such that
\[\|f\|_{B_{q,q}^{(\frac1q-\frac12)d}(\R^d;E)} \leq C'\|f|_{B_{2r}}\|_{B_{q,q}^{(\frac1q-\frac12)d}(B_{2r};E)}.\]
Choosing $n$ so large that $B_{2r} \subseteq D_n$, we may conclude that
\[
\|f\|_{B_{q,q}^{(\frac1q-\frac12)d}(\R^d;E)} \leq
C'\|f_n\|_{B_{q,q}^{(\frac1q-\frac12)d}(D_n;E)}\leq C'C \|f\|_{\g(\R^d;E)}.
\]
Since $C^{\infty}_c(\R^d;E)$ is dense in $\g(\R^d;E)$ the result follows from
Theorem \ref{thm:main}.
\end{proof}

It is an interesting fact that at least in dimension $d=1$, the ``if part" of
Theorem \ref{thm:embed} (1) can be improved as follows.

\begin{theorem}\label{thm:counterpartBesov}
If $p\in [1,2)$ is such that we have a continuous embedding $$B^{\frac{1}{p}-\frac12}_{p,1}((0,1);E)\embed
\g(L^2(0,1), E),$$ then $E$ has type $p$.
\end{theorem}

\begin{proof}
We may assume that $p\in (1,2)$.

First, for $s>0$ we introduce an equivalent norm on $B_{p,q}^s(\R;E)$ which
does not involve the Fourier transform and can be handled quite easily from
the computational point of view.

For $h\in\R$ and a function $f:\R\to E$ we define the function $T(h)f:\R\to E$
as the translate of $f$ over $h$, i.e.\
$$(T (h) f)(t) := f(t+h).$$
For $f\in L^p(\R;E)$ and
$t>0$ let
\[\varrho_p(f,t) := \sup_{|h|\le t} \|T(h) f - f\|_{L^p(\R;E)}.\]
Then
$$
\|f\|_{B^s_{p, q}(\R;E)}^* := \|f\|_{L^p(\R;E)}
+ \Big(\int_0^1 \big(t^{-s}\varrho_p(f,t)\big)^q\, \frac{dt}{t}\Big)^\frac1q
$$
(with the obvious modification for $q=\infty$) defines an equivalent norm on
$B^s_{p, q}(\R;E)$ (see \cite[Proposition 3.1]{PeWo} or \cite[Theorem
4.3.3]{Schm}).

With these preliminaries out of the way we turn to the proof of the theorem.
Since every Banach space has type 1 we may assume that $p\in (1,2)$. Let
$n\geq 1$ and $x_0, \ldots, x_{n-1}\in E$ be arbitrary and fixed. For $j=0,
\ldots, 2n-1$, let $t_j = \frac{j}{2n}$. Define $f:\R\to E$ as
\[f = \sum_{k=0}^{n-1} \one_{(t_{2k}, t_{2k+1}]} x_k.\]
Then $\|f\|_{L^p(\R;E)} = (2n)^{-\frac1p} \big(\sum_{k=0}^{n-1}
\|x_k\|^p\big)^{\frac1p}$. Let $0<t<(2n)^{-1}$ and take $0<|h|\le t$. If
$h>0$, then
\[T(h) f - f = 
\sum_{k=0}^{n-1} (\one_{(t_{2k}-h, t_{2k}]} - \one_{(t_{2k+1}-h, t_{2k+1}]})x_k.\] 
If $h<0$, then
\[T(h) f - f = \sum_{k=0}^{n-1} (-\one_{(t_{2k}, t_{2k}+h]}+ \one_{(t_{2k+1}, t_{2k+1}+h]})x_k.
\] In both cases we find that
\[
\|T(h) f - f\|_{L^p(\R;E)}^p \leq  2|h| \sum_{k=0}^{n-1} \|x_k\|^p\leq 2 t
\sum_{k=0}^{n-1} \|x_k\|^p.\] This shows that $\varrho_p(f, t) \leq
2^{\frac1p} t^{\frac1p} \big(\sum_{k=0}^{n-1} \|x_k\|^p\big)^{\frac1p}$ for
all $0<t<(2n)^{-1}$. It follows that
\[\begin{aligned}
\int_0^{(2n)^{-1}}  t^{-\frac1p + \frac12}\varrho_p(f, t) \frac{dt}{t} & \leq
2^{\frac1p} \Big(\sum_{k=0}^{n-1} \|x_k\|^p\Big)^{\frac1p} \int_0^{(2n)^{-1}}
t^{\frac{1}{2}} \frac{dt}{t}
\\ &= 2^{\frac1p+1} (2n)^{-\frac{1}{2}}
\Big(\sum_{k=0}^{n-1} \|x_k\|^p\Big)^{\frac{1}{p}}.
\end{aligned}
\] If $t>(2n)^{-1}$, then $\varrho_p(f, t)\leq 2\n f\n_p = 2 (2n)^{-\frac{1}{p}}
\big(\sum_{k=0}^{n-1} \|x_k\|^p\big)^{\frac1p}$. It follows that
\[\begin{aligned}
\int_{(2n)^{-1}}^1 t^{-\frac1p +\frac12}\varrho_p(f, t) \frac{dt}{t} & \leq 2
(2n)^{-\frac{1}{p}} \Big(\sum_{k=0}^{n-1} \|x_k\|^p\Big)^{\frac1p}
\int_{(2n)^{-1}}^1 t^{-\frac1p +\frac12} \frac{dt}{t} \\ & = 2
(2n)^{-\frac{1}{p}} \Big(\sum_{k=0}^{n-1} \|x_k\|^p\Big)^{\frac1p}
\frac{1}{\frac1p-\frac12} ((2n)^{\frac1p -\frac12} - 1)
\\ & \leq 2 (2n)^{-\frac{1}{2}}
\frac{1}{\frac1p-\frac12} \Big(\sum_{k=0}^{n-1} \|x_k\|^p\Big)^{\frac1p}.
\end{aligned}\]
It follows that $f\in B^{\frac1p - \frac12}_{p, 1}(\R;E)$ and by restricting
to $(0,1)$ we obtain
\[\|f\|_{B^{\frac1p - \frac12}_{p, 1}((0,1);E)} \le \|f\|_{B^{\frac1p - \frac12}_{p, 1}(\R;E)}\leq C_p (2n)^{-\frac{1}{2}} \Big(\sum_{k=0}^{n-1} \|x_k\|^p\Big)^{\frac1p},\]
where $C_p$ depends only on $p$.
On the other hand,
\[\|I_{f}\|_{\g(L^2(0,1), E)} = (2n)^{-\frac{1}{2}} \Big\|\sum_{k=0}^{n-1}
\g_k x_k\Big\|_{L^2(\O;E)}.\]
From the boundedness of the embedding $I: B_{p,1}^{\frac1p-\frac12}((0,1);E)\embed \g(L^2(0,1), E)$ we conclude that
\[(2n)^{-\frac{1}{2}} \Big\|\sum_{k=0}^{n-1} \g_k x_k\Big\|_{L^2(\O;E)} \leq  C_p (2n)^{-\frac{1}{2}}\n I\n \Big(\sum_{k=0}^{n-1} \|x_k\|^p\Big)^{\frac1p}.\]
Hence $E$ has type $p$, with Gaussian type $p$ constant of at most $C_p\n I\n$.
\end{proof}


Returning to Theorem \ref{thm:embed}, we note the following consequence:

\begin{corollary}\label{cor:holder}
Let $D\subseteq \R^d$ be a nonempty bounded open domain
with smooth boundary. Let $p\in [1,2]$ and $\a,\b\in\R$ satisfy $\a> \b +(\frac{1}{p} - \frac{1}{2})d\ge 0$.
If $E$ has type $p$, we have a continuous embedding
$$ C^{\a}(\overline D;E)\embed \g(H^{\b,2}(D),E).$$
\end{corollary}
\begin{proof}
For $\a>\g>\b +(\frac{1}{p} - \frac{1}{2})d\ge 0$ we have, cf. \cite{Am2},
$$ C^{\a}(\overline D;E)\embed B_{\infty,\infty}^\g(D;E)\embed
B_{p,p}^{\b+(\frac{1}{p} - \frac{1}{2})d}(D;E).$$ The result now follows from Theorem \ref{thm:embed}.
\end{proof}

For dimension $d=1$ we have the following converse:

\begin{theorem}\label{thm:counterpart}
Let $E$ be a  Banach space, and let  $p\in (1,2)$ and $\a\in (0,\frac1p -
\frac12)$. If $C^{\a}([0,1];E)\embed \g((0,1);E)$,
then $E$ has type $p$.
\end{theorem}

In particular this shows that in the spaces
$E=l^p$ and $E=L^p(0,1)$, with $p\in
[1, 2)$, for every $ \alpha \in (0,\frac1p - \frac12)$
there exist $\alpha$-H\"older continuous functions
which do not belong to $\g(L^2(0,1), E)$.
Indeed, for such $\a$ we can find $p<p'<2$ such that $\alpha\in (0, \frac1{p'} - \frac12)$, but both $l^p$ and $L^p(0,1)$ fail type $p'$.
A similar result holds for $E=c_0$ and $E=C([0,1])$  and $\a\in (0,\frac12)$.
This improves the examples in \cite{RS}, where only measurable functions are considered.

\begin{proof}
Assume for a contradiction that $E$ is not of type $p$. We will show that this
leads to a contradiction. By the Maurey-Pisier theorem (see \cite{MP}), $l^p$
is finitely representable in $E$. Fix an integer $n$ and let $T:l^p_n \to E$
be such that for all $x\in l^p_n$
\[\|x\|_{l^p_n}\leq \|T x\|\leq 2\|x\|_{l^p_n}.\]
Choose $1<r<(\frac12p+\alpha p)^{-1}$. Let $c= \sum_{i\geq 1} i^{-r}$ and let
$t_0 = 0$, $t_k = c^{-1} \sum_{i=1}^k i^{-r}$ for $k\geq 1$. Let
$(e_k)_{k=1}^n$ be the standard basis of $l^p_n$ and define $g_n:[0,1]\to
l^p_n$ as
\[g_n(t) = \left\{
\begin{array}{ll}
    \Big(1-\frac{|2t-t_k - t_{k-1}|}{t_k - t_{k-1}}\Big) e_k , & \text{if $t\in (t_{k-1}, t_k]$ for $1\leq k\leq n$,} \\
    0, & \text{otherwise.} \\
\end{array}
\right.\] We claim that $g_n$ is H\"older continuous of exponent $\a$ and
\[\begin{aligned}
\|g_n\|_{C^\alpha([0,1];l^p_n)} & = \sup_{t\in [0,1]}\|g_n(t)\|_{l^p_n} +
\sup_{0\leq s<t\leq 1} \frac{\|g(t)-g(s)\|_{l^p_n}}{|t-s|^{\alpha}}\\ & \leq 1
+ 4 (t_n - t_{n-1})^{-\alpha} = 1+ 4 c^{\alpha} n^{r \alpha}.
\end{aligned}
\]
To show this we consider several cases. First of all $\|g_n(t)\|_{l^p_n}\leq
1$ for all $t\in [0,1]$. If $t,s\in [t_{k-1}, t_{k}]$ for some $1\leq k\leq
n$,
\[\begin{aligned}
\|g_n(t)-g_n(s)\|_{l^p_n} &= \Big|\frac{|2t-t_k - t_{k-1}|}{t_k - t_{k-1}} -
\frac{|2s-t_k - t_{k-1}|}{t_k - t_{k-1}}\Big|\\ & \leq
\frac{2|t-s|}{t_k-t_{k-1}}\leq \frac{2|t-s|^\alpha}{|t_k-t_{k-1}|^\alpha}\leq
\frac{2|t-s|^\alpha}{|t_n-t_{n-1}|^\alpha}.
\end{aligned}\]
If $s\in (t_{k-1}, t_k]$ and $t\in (t_k, t_{k+1}]$ for some $1\leq k\leq n-1$,
then by the above estimate and the concavity of $x\mapsto x^{\alpha}$,
\[\begin{aligned}
\|g_n(t)-g_n(s)\|_{l^p_n}&\leq \|g_n(t)-g_n(t_k)\|_{l^p_n} +
\|g_n(t_k)-g_n(s)\|_{l^p_n}
\\ & \leq \frac{2|t-t_k|^\alpha}{|t_n-t_{n-1}|^\alpha} + \frac{2|t_k-s|^\alpha}{|t_n-t_{n-1}|^\alpha}\leq \frac{ 2^{2-\alpha} |t-s|^\alpha}{|t_n-t_{n-1}|^\alpha}.
\end{aligned}\]
If $s\in (t_{l-1}, t_l]$ and $t\in (t_{k-1}, t_k]$ for $l+2\leq k\leq n$ then
\[\|g_n(t)-g_n(s)\|_{l^p_n}\leq 2 \leq 2 (t_{k-1}-t_l)^{\alpha} (t_n-t_{n-1})^{-\alpha}\leq 2 (t-s)^{\alpha} (t_n-t_{n-1})^{-\alpha}.\]
For the other cases the estimate is obvious and we proved the claim. We have
$g_n\in \g(L^2(0,1),l^p_n)$ and a standard square function estimate (cf.
\cite[Example 7.3]{NW}) gives
\[\begin{aligned}
\|I_{g_n}\|^p_{\g(L^2(0,1),l^p_n)} & \ge K_{p}^p \sum_{k=1}^n
\Big(\int_{t_{k-1}}^{t_k} \Big(1-\frac{|2t-t_k - t_{k-1}|}{t_k -
t_{k-1}}\Big)^2 \, dt\Big)^{\frac{p}{2}}
\\ & = K_{p}^p\sum_{k=1}^n
\Big(\frac{t_k-t_{k-1}}{2}\int_{-1}^{1} (1-|s|)^2 \, dt\Big)^{\frac{p}{2}}
\\ & \ge  3^{-\frac{p}{2}} c^{-\frac{p}{2}}K_{p}^p \frac{2}{2-pr}( (n+1)^{-\frac{p
r}{2}+1} - 1 ),
\end{aligned}\]
where $K_p$ is a constant depending only on $p$. Define $f_n:[0,1] \to E$ as
$f_n := T g_n$. Then $f_n$ is $\alpha$-H\"older continuous and $I_{f_n}\in
\g(L^2(0,1), E)$ with
\[\|f_n\|_{C^\alpha([0,1], E)}\leq
2\|g_n\|_{C^\alpha([0,1], l^p_n)}\leq 2(1+ 4c^{\alpha} n^{r \alpha}).\] and
\[\|I_{f_n}\|_{\g(L^2(0,1), E)}^p \geq  \|I_{g_n}\|_{\g(L^2(0,1), l^p_n)}^p \ge 3^{-\frac{p}{2}} c^{-\frac{p}{2}}
K_{p}^p\frac{2}{2-pr}( (n+1)^{-\frac{p r}{2}+1}
- 1 ).\]
Since the inclusion operator $I : C^{\a}([0,1];E)\to \g(L^2(0,1),E)$ is bounded we conclude that
\[3^{-\frac{1}{2}} c^{-\frac{1}{2}} K_{p}\frac{2^{\frac1p}}{(2-pr)^{\frac1p}}( (n+1)^{-\frac{p r}{2}+1}
- 1 )^{\frac1p} \le 2(1+ 4 c^{\alpha} n^{r \alpha}).\] Since we may take
$n$ arbitrary large, this implies $-\frac{r}{2}+\frac1p\leq r\alpha$, so $r\geq
(\alpha p+ \frac{p}{2})^{-1}$. But this contradicts the choice of $r$, and the proof is complete.
\end{proof}

\providecommand{\bysame}{\leavevmode\hbox to3em{\hrulefill}\thinspace}

\end{document}